\documentclass[11pt]{article}
\usepackage{amssymb,amsmath,amsfonts}
\usepackage{color}

\begin{document}

\author{Constantin C\u alin and Mircea Crasmareanu}
\title{From the Eisenhart problem to Ricci solitons in $f$-Kenmotsu manifolds}
\date{Dedicated to the memory of Neculai Papaghiuc 1947-2008}
\maketitle

\begin{abstract}
The Eisenhart problem of finding parallel tensors is solved for the symmetric case in the regular $f$-Kenmotsu framework. On this way, the Olszack-Rosca example of Einstein manifolds provided by $f$-Kenmotsu manifolds via locally symmetric Ricci tensors is recovered as well as a case of Killing vector fields. Some other classes of Einstein-Kenmotsu manifolds are presented. Our result is interpreted in terms of Ricci solitons and special quadratic first integrals.
\end{abstract}



\medskip

\section*{Introduction}

In 1923, Eisenhart \cite{l:e} proved that if a positive definite Riemannian ma\-ni\-fold $(M,g)$ admits a second order parallel symmetric covariant tensor other that a constant multiple of the metric tensor, then it is reducible. In 1926, Levy \cite{h:l} proved that a second order parallel symmetric non-degenerated tensor $\alpha $ in a space form is proportional to the metric tensor. Let us point out that this question can be considered as dual to the the problem of finding linear connections making parallel a given tensor field, problem which was considered by Wong in \cite{y:w}. Also, the former question implies topological restrictions namely if the (pseudo) Riemannian manifold $M$ admits a parallel symmetric $(0,2)$ tensor field then $M$ is locally the direct product of a number of (pseudo) Riemannian manifolds, \cite{h:w} (cited by \cite{g:z}). Another situation where the parallelism of $\alpha $ is involved appears in the theory of {\it totally geodesic maps}, namely, as is point out in \cite[p. 114]{o:c}, $\nabla \alpha =0$ is equivalent with the fact that $1:(M,g)\rightarrow (M,\alpha )$ is a totally geodesic map.

While both Eisenhart and Levy work locally, Ramesh Sharma gives in \cite{r:s1} a global approach based on Ricci identities. In addition to space-forms, Sharma considered this {\it Eisenhart problem} in contact geometry \cite{r:s2}-\cite{r:s4}, for example for $K$-contact manifolds in \cite{r:s3}. Since then, several other studies appear in various contact manifolds: nearly Sasakian \cite{t:m}, (para) $P$-Sasakian \cite{t:d1}, \cite{d:e} and \cite{z:l}, $\alpha $-Sasakian \cite{l:d2}. Another framework was that of quasi-constant curvature in \cite{x:j}. Also, contact metrics with nonvanishing $\xi $-sectional curvature are studied in \cite{g:s}.

Returning to contact geometry, an important class of manifolds are introduced by Kenmotsu in \cite{k:k} and generalized by Olszack and Rosca in \cite{o:r}. In the last time, there is an increasing flow of papers in this direction e.g. that of our professor N. Papaghiuc \cite{n:p1}-\cite{n:p2} to which we dedicate this short note. Motivated by this fact we studied the case of $f$-Kenmotsu manifolds satisfying a special condition called by us {\it regular} and show that a symmetric parallel tensor field of second order must be a constant multiple of the Riemannian metric. There are three remarks regarding our result:\\
i) it is in agreement with what happens in all previously recalled contact \\ geometries for the symmetric case, \\
ii) it is obtained in the same manner {\bf as} originated in Sharma's paper \cite{r:s1}, \\
iii) yields a class of Einstein manifolds already indicated by Olszack and Rosca but with a more complicated proof. \\
Let us point out also that the anti-symmetric case appears without proof in \cite{v:m}.

Our main result is connected with the recent theory of Ricci solitons, a subject included in the Hamilton-Perelman approach (and proof)  of Poincar\'e Conjecture. Ricci solitons in contact geometry were first studied by Ramesh Sharma in \cite{g:sc} and \cite{r:s5}; also the preprint \cite{mm:t} is available to arxiv. In these papers the $K$-contact and $(k, \mu )$-contact (including Sasakian) cases are treated; then our treatment for the Kenmotsu variant of almost contact geometry seems to be new.

Our work is structured as follows. The first section is a very brief review of Kenmotsu geometry and Ricci solitons. The next section is devoted to the (symmetric case of) Eisenhart problem in a $f$-Kenmotsu manifold and several situations yielding Einstein manifolds are derived. Also, the relationship with the Ricci solitons is pointed out. The last section offers a dynamical picture of the subject via Killing vector fields and quadratic first integrals of a special type.

{\bf Acknowledgement} Special thanks are offered to Gheorghe Pitis for some useful remarks as well as sending us his book \cite{gh:p}, a source of several references. Also, we are very indebted to Marian-Ioan Munteanu and the referees who pointed out major improvements.

\section{$f$-Kenmotsu manifolds. Ricci solitons}

Let $M$ be a real $2n+1$-dimensional differentiable manifold endowed with an almost
contact metric structure $(\varphi, \xi, \eta, g)$:
$$
\begin{array}{l}
(a) \quad \varphi ^{2} = -I + \eta \otimes \xi,\ \ (b) \quad \eta(\xi) = 1,\ \ (c) \quad \eta \circ \varphi = 0,\\ (d) \quad \varphi (\xi) = 0, \ \
(e) \quad \eta(X) = g(X, \xi),\\  (f) \quad g(\varphi X, \varphi Y)=g(X, Y)-\eta (X)\eta(Y),
\end{array} \eqno(1.1)
$$
for any vector fields $X, Y\in {\cal X}(M)$ where $I$ is the identity of the tangent
bundle $TM$, $\varphi $ is a tensor field of $(1,1)$-type, $\eta $ is a 1-form,
$\xi $ is a vector field and $g$ is a metric tensor field.
Throughout the paper all objects are differentiable of class $C^{\infty }$.

\medskip

We say that $(M,\varphi ,\xi ,\eta ,g)$ is an $f$-{\it Kenmotsu manifold} if the Levi-Civita connection of $g$ satisfy \cite{v:m}:
$$
(\nabla _X\varphi )(Y)=f(g(\varphi X,Y)\xi -\varphi (X)\eta (Y)), \eqno(1.2)
$$
where $f\in C^{\infty }(M)$ is strictly positive and $df\wedge \eta =0$ holds. A $f=constant\equiv \beta >0$ is called $\beta $-{\it Kenmotsu manifold} with the particular case $f\equiv 1$-Kenmotsu manifold which is a usual {\it Kenmotsu manifold} \cite{k:k}.

\medskip

In a general $f$-Kenmotsu manifold we have, \cite{o:r}:
$$
\nabla _X\xi =f(X-\eta (X)\xi), \eqno(1.3)
$$
and the curvature tensor field:
$$
R(X,Y)\xi =f^2(\eta (X)Y-\eta (Y)X)+Y(f)\varphi ^2X-X(f)\varphi ^2Y \eqno(1.4)
$$
while the Ricci curvature and Ricci tensor are, \cite{k:pt}:
$$
S(\xi ,\xi ) = -2n(f^2+\xi (f)) \eqno(1.5)
$$
$$
Q(\xi )=-2nf^2\xi -\xi (f)\xi-(2n-1)grad f. \eqno(1.6)
$$

In the last part of this section we recall the notion of Ricci solitons according to \cite[p. 139]{r:s5}. On the manifold $M$, a {\it Ricci soliton} is a triple $(g,V,\lambda )$ with $g$ a Riemannian metric, $V$ a vector field and $\lambda $ a real scalar such that:
$$
{\cal L}_Vg+2S+2\lambda g=0. \eqno(1.7)
$$
The Ricci soliton is said to be {\it shrinking, steady} or {\it expanding} according as $\lambda $ is negative, zero or positive.

\section{Parallel symmetric second order tensors and \\ Ricci solitons in $f$-Kenmotsu manifolds}

Fix $\alpha $ a symmetric tensor field of $(0,2)$-type which we suppose to be parallel with respect to $\nabla $ i.e. $\nabla \alpha =0$. Applying the Ricci identity \\ $\nabla ^2\alpha (X,Y;Z,W)-\nabla ^2(X,Y;W,Z)=0$ we get the relation $(1.1)$ of \cite[p. 787]{r:s1}:
$$
\alpha (R(X,Y)Z,W) +\alpha (Z,R(X,Y)W)=0, \eqno(2.1)
$$
which is fundamental in all papers treating this subject. Replacing $Z=W=\xi $ and using $(1.4)$ it results:
$$
f^2[\eta (X)\alpha (Y,\xi )-\eta (Y)\alpha (X,\xi)]+Y(f)\alpha
(\varphi ^2X, \xi)-X(f)\alpha (\varphi ^2Y, \xi)=0, \eqno(2.2)
$$
by the symmetry of $\alpha $. With $X=\xi $ we derive:
$$
[f^2+\xi (f)][\alpha (Y,\xi )-\eta (Y)\alpha (\xi , \xi)]=0
$$
and supposing $f^2+\xi (f)\neq 0$ it results:
$$
\alpha (Y,\xi )=\eta (Y)\alpha (\xi , \xi). \eqno(2.3)
$$
Let us call {\it regular} $f$-{\it Kenmotsu manifold} a $f$-Kenmotsu manifold with $f^2+\xi (f)\neq 0$ and remark that $\beta $-Kenmotsu manifolds are regular.

\smallskip

Differentiating the last equation covariantly with respect to $X$ we have:
$$
\alpha (\nabla _XY, \xi )+f[\alpha (X,Y)-\eta (X)\eta (Y)\alpha (\xi ,\xi )]=X(\eta (Y))\alpha (\xi , \xi), \eqno(2.4)
$$
which means via $(2.3)$ with $Y \rightarrow \nabla _XY$:
$$
f[\alpha (X,Y)-\eta (X)\eta (Y)\alpha (\xi ,\xi )]=[X(g(Y,\xi ))-g(\nabla _XY, \xi )]\alpha (\xi ,\xi )=
$$
$$
=g(Y,\nabla _X\xi )\alpha (\xi ,\xi )=
f[g(X,Y)-\eta (X)\eta (Y)]\alpha (\xi ,\xi ). \eqno(2.5)
$$
From the positiveness of $f$ we deduce that:
$$
\alpha (X,Y)=\alpha (\xi ,\xi )g(X,Y) \eqno(2.6)
$$
which together with the standard fact that the parallelism of $\alpha $ implies the constance of $\alpha (\xi ,\xi )$ via $(2.3)$ yields:

\medskip

{\bf Theorem} {\it A symmetric parallel second order covariant tensor in a regular $f$-Kenmotsu manifold is a constant multiple of the metric tensor. In other words, a regular $f$-Kenmotsu metric is irreducible which means that the tangent bundle does not admits a decomposition $TM=E_1\oplus E_2$ parallel with respect of the Levi-Civita connection of $g$}.

\medskip

{\bf Corollary 1} {\it A locally Ricci symmetric} $(\nabla S\equiv 0)$ {\it regular $f$-Kenmotsu manifold is an Einstein manifold}.

\medskip

{\bf Remarks} 1) The particular case of dimension three and $\beta$-Kenmotsu of our theorem appears in Theorem 3.1 from \cite[p. 2689]{d:m}. The above corollary has been proved by Olszack and Rosca in another way. \\
2) In \cite{b:tdt} it is shown the equivalence of the following statements for an Kenmotsu manifold:\\
i) is Einstein, \\
ii) is locally Ricci symmetric, \\
iii) is Ricci semi-symmetric i.e. $R\cdot S=0$ where:
$$
(R(X, Y)\cdot S)(X_1, X_2) = -S(R(X, Y)X_1, X_2)-S(X_1, R(X, Y)X_2).
$$
The same implication iii) $\rightarrow $ i) for Kenmotsu manifolds is Theorem 1 from \cite[p. 438]{j:dp}. But we have the implication iii) $\rightarrow $ i) in the more general framework of regular $f$-Kenmotsu manifols since $R\cdot S=0$ means exactly $(2.1)$ with $\alpha $ replaced by $S$. Every semisymmetric manifold, i. e. $R\cdot R=0$, is Ricci-semisymmetric but the converse statement is not true. In conclusion:

\medskip

{\bf Proposition 1} {\it A Ricci-semisymmetric, particularly semisymmetric, re\-gu\-lar} $f$-{\it Kenmotsu manifold is Einstein.}

\medskip

Another class of spaces related to the Ricci tensor was introduced in \cite{h:q}; namely a Riemannian manifold is a {\it special weakly Ricci symmetric space} if there exists a 1-form $\rho $ such that:
$$
(\nabla _XS)(Y,Z)=2\rho (X)S(Y,Z)+\rho (Y)S(Z,X)+\rho (Z)S(X,Y). \eqno(2.7)
$$
The same condition was sometimes called {\it generalized pseudo-Ricci sym\-metric manifold} (\cite{j:s}) or simply {\it pseudo-Ricci symmetric manifold} (\cite{c:k}). Making $X=Y=Z=\xi $ it results:
$$
\xi (S(\xi ,\xi ))=4\rho (\xi )S(\xi ,\xi) \eqno(2.8)
$$
and then for a $\beta $-Kenmotsu manifold we get $\rho (\xi )=0$. Returning to $(2.7)$ with $Y=Z=\xi $ will results $\rho (X)=0$ for every vector field $X$ and then we have a generalization of Theorem 3.3. from \cite[p. 96]{a:go}:

\medskip

{\bf Proposition 2} {\it A} $\beta $-{\it Kenmotsu manifold which is special weakly Ricci symmetric is an Einstein space.}

\medskip

We close this section with applications of our Theorem to Ricci solitons:

\medskip

{\bf Corollary 2} {\it Suppose that on a regular $f$-Kenmotsu manifold the} $(0,2)$-{\it type field} ${\cal L}_Vg+2S$ {\it is parallel where} $V$ {\it is a given vector field. Then} $(g, V)$ {\it yield a Ricci soliton. In particular, if the given regular} $f$-{\it Kenmotsu manifold is Ricci-semisymmetric or semisymmetric with} ${\cal L}_Vg$ {\it parallel, we have the same conclusion.}

\medskip

Naturally, two situations appear regarding the vector field $V$: $V\in span{\xi }$ and $V\bot \xi$ but the second class seems far too complex to analyse in practice. For this reason it is appropriate to investigate only the case $V=\xi $.

We are interested in expressions for ${\cal L}_{\xi }g+2S$. A straightforward computation gives:
$$
{\cal L}_{\xi }g(X,Y)=2f(g(X,Y)-\eta (X)\eta (Y))=2fg(\varphi X, \varphi Y). \eqno(2.9)
$$

A general expression of $S$ is known by us only for the the 3-dimensional case and $\eta $-Einstein Kenmotsu manifolds. Let us treat these situations in the following:

\medskip

I) \cite[p. 251]{d:t}:
$$
S(X,Y)=\left(\frac{r}{2}+\xi (f)+f^2\right)g(X,Y)-
$$
$$
-\left(\frac{r}{2}+\xi (f)+3f^2\right)\eta (X)\eta (Y)-Y(f)\eta (X)-X(f)\eta (Y) \eqno(2.10)
$$
where $r$ is the scalar curvature. Then, for a 3-dimensional $f$-Kenmotsu manifold we get:
$$
\alpha:=({\cal L}_{\xi }g+2S)(X,Y)=(r+2\xi (f)+2f+2f^2)g(X,Y)-
$$
$$
-(r+2\xi (f)+2f+6f^2)\eta (X)\eta (Y)-2Y(f)\eta (X)-2X(f)\eta (Y) \eqno(2.11)
$$
while, for $\beta $-Kenmotsu:
$$
\alpha (X,Y)=(r+2\beta +2\beta ^2)g(\varphi X, \varphi Y)-4\beta ^2\eta (X)\eta (Y), \eqno(2.12)
$$
$$
(\nabla _Z\alpha )(X,Y)=Z(r)g(\varphi X, \varphi Y)-
$$
$$
-\beta (r+2\beta +6\beta ^2)[\eta (X)g(\varphi Y, \varphi Z)+\eta (Y)g(\varphi X, \varphi Z)]. \eqno(2.13)
$$
{\bf Substituting} $Z=\xi , X=Y\in (span{\xi })^{\bot }$ respectively $X=Y=Z\in (span{\xi })^{\bot }$ in $(2.13)$ we derive that $r$ is a constant, provided $\alpha $ is parallel. Thus, we can state the following:

\medskip

{\bf Proposition 3} {\it A 3-dimensional} ${\beta }$-{\it Kenmotsu Ricci soliton} $(g, \xi , \lambda )$ {\it is expanding and with constant scalar curvature}.

\medskip

{\bf Proof} $\lambda =-\frac{1}{2}\alpha (\xi, \xi)=2\beta ^2$. \quad $\Box $

\medskip

At this point we remark that the Ricci solitons of almost contact geometry studied in \cite{r:s5} and \cite{mm:t} in relationship with the Sasakian case are shrinking and this observation is in accordance with the diagram of Chinea from \cite{d:c} that Sasakian and Kenmotsu are opposite sides of the trans-Sasakian moon. Also, the expanding character may be considered as a manifestation of the fact that a $\beta $-Kenmotsu manifold can not be compact.

\medskip

II) Recall that the metric $g$ is called $\eta $-{\it Einstein} if there exists two real functions $a, b$ such that the Ricci tensor of $g$ is:
$$
S=ag+b\eta \otimes \eta .
$$
For an $\eta $-Einstein Kenmotsu manifold we have, \cite[p. 441]{j:dp}:
$$
S(X,Y)=\left(\frac{r}{2n}+1\right)g(X,Y)-\left(\frac{r}{2n}+2n+1\right)\eta (X)\eta (Y) \eqno(2.14)
$$
and then:
$$
\alpha (X,Y)=\left(\frac{r}{n}+4\right)g(X,Y)-\left(\frac{r}{n}+4+4n\right)\eta (X)\eta (Y) \eqno(2.15)
$$
$$
(\nabla _Z\alpha )(X,Y)=\frac{1}{n}Z(r)g(\varphi X, \varphi Y)-
$$
$$
-\left(\frac{r}{n}+4n+4\right)[\eta (Y)g(\varphi X, \varphi Z)+\eta (X)g(\varphi Y, \varphi Z)]. \eqno(2.16)
$$

\medskip

{\bf Proposition 4} {\it An} ${\eta }$-{\it Einstein Kenmotsu Ricci soliton} $(g, \xi , \lambda )$ {\it is ex\-pan\-ding and with constant scalar curvature, thus Einstein}.

\medskip

{\bf Proof} $\lambda =-\frac{1}{2}\alpha (\xi, \xi)=2n$. The same computation as in Proposition 3 implies constant scalar curvature. \quad $\Box $

\section{The dynamical point of view}

We begin this section with a straightforward consequence of the main Theorem, which also appears in the Olzack-Rosca paper, and is related to the last part of Corollary 2:

\medskip

{\bf Corollary 3} {\it An affine Killing vector field in a $\beta $-Kenmotsu manifold is Killing. As consequence, that scalar provided by the Ricci soliton} $(g, V)$ {\it of a Ricci-semisymmetric} $\beta $-{\it Kenmotsu manifold is} $\lambda =-S(V,V)$.

\medskip

{\bf Proof} (inspired by \cite[p. 504]{g:s}) Fix $X\in {\cal X}(M)$ an affine Killing vector field: $\nabla {\cal L}_Xg=0$. From Theorem it results that $X$ is {\it conformal Killing} i.e. ${\cal L}_Xg = cg$; more precisely $X$ is {\it homothetic} since $c$ is a constant. Lie differentiating the identity
$(1.5)$ along $X$ and using ${\cal L}_XS = 0$ (since $X$ is homothetic) and equation $(1.6)$
we get $g({\cal L}_X\xi ,\xi ) = 0$. Hence $c =({\cal L}_Xg)(\xi ,\xi )=-2g({\cal L}_X\xi ,\xi )= 0$. Thus $X$ is Killing. \quad $\Box $

\medskip

Let us present another dynamical picture of our results. Let $(M, \nabla)$ be a $m$-dimensional manifold endowed with a symmetric linear connection. A {\it quadratic first integral} (QFI on short) for the geodesics of $\nabla $ is defined by ${\cal F}=a_{ij}\frac{dx^i}{dt}\frac{dx^j}{dt}$ with a symmetric 2-tensor field $a=(a_{ij})$ satisfying the {\it Killing-type equations}:
$$
a_{ij:k}+a_{jk:i}+a_{ki:j}=0, \eqno(3.1)
$$
where, as usual, the double dot means the covariant derivative with respect to $\nabla $.

\medskip

The QFI defined by $a$ is called {\it special} (SQFI) if
$a_{ij:k}=0$ and the maximum number of linearly independent SQFI a pair $(M, \nabla)$ can admit
is $\frac{m(m+1)}{2}$; a flat space will admit this number.
In \cite[p. 117]{l:k} it is shown that a non-flat Riemannian
manifold may admit as many as $M_S(m)=1+\frac{(m-2)(m-1)}{2}$
linearly independent SQFI. Therefore, for an almost contact manifold
($m=2n+1$) the maximum number of SQFI is $M_S(2n+1)=1+n(2n-1)>1$.

\medskip

Our main result implies that for a regular $f$-Kenmotsu manifold the number of SQFI is exactly 1 and the only SQFI is the {\it kinetic energy} ${\cal F}=g_{ij}\frac{dx^i}{dt}\frac{dx^j}{dt}$. So:

\medskip

{\bf Proposition 5} {\it There exist almost contact manifolds which does not admit $M_S(2n+1)$ SQFI}.

\medskip

It remains as an open problem to find examples of almost contact metrics with exactly $M_S(2n+1)$ SQFI.

\medskip

\noindent Department of Mathematics, \newline
Technical University "Gh.Asachi"\newline
Ia\c si, 700049\newline
Romania\newline
e-mail: c0nstc@yahoo.com

\medskip

\noindent Faculty of Mathematics
\newline University "Al. I.Cuza" \newline
Ia\c si, 700506 \newline Romania \newline e-mail: mcrasm@uaic.ro

\smallskip

\noindent http://www.math.uaic.ro/$\sim$mcrasm

\end{document}